\documentclass[12pt]{amsart}
\usepackage{url}
\usepackage{amsmath,amsxtra,amssymb,latexsym,epsfig,amscd,amsthm,fancybox,epsfig}
\usepackage[mathscr]{eucal}
\usepackage{graphicx}
\usepackage{multicol,xcolor}
\usepackage{epsfig} 
\usepackage{epstopdf}
\usepackage{cases}
\usepackage[subfigure]{graphfig}

\usepackage{color}
\usepackage{hyperref}
\usepackage{bigints}
\usepackage{enumitem}
\usepackage{bm}
\usepackage{bbm}

\usepackage{natbib}
\usepackage{hyperref}
\hypersetup{
    colorlinks=true,
    linkcolor=blue,
    filecolor=magenta,
    urlcolor=cyan,
    citecolor=blue,
}

\usepackage{thmtools,thm-restate}

\newtheorem{thm}{Theorem}[section]

\newtheorem{rmk}{Remark}[section]

\newtheorem{lem}{Lemma}[section]

\theoremstyle{definition}

\theoremstyle{remark}

\numberwithin{equation}{section}

\newcommand{\E}{\mathbb{E}}

\newcommand{\BX}{\mathbf{X}}

\newcommand{\ba}{\mathbf{a}}
\newcommand{\bb}{\mathbf{b}}
\newcommand{\bc}{\mathbf{c}}

\numberwithin{equation}{section}

\newcommand{\bed}{\begin{displaymath}}
\newcommand{\eed}{\end{displaymath}}
\newcommand{\bea}{\bed\begin{array}{rl}}
\newcommand{\eea}{\end{array}\eed}

\newcommand{\barray}{\begin{array}{ll}}
\newcommand{\earray}{\end{array}}

\newcommand{\1}{\boldsymbol{1}}

\definecolor{OliveGreen}{cmyk}{0.64,0,0.95,0.40}

\def\bar{\overline}

\def\a.s{\text{\;a.s.\;}}

\newcommand{\bmu}{{\boldsymbol\mu}}
\newcommand{\BSigma}{{\boldsymbol\Sigma}}
\newcommand{\bee}{\mathbf{e}}
\newcommand{\bk}{\mathbf{k}}
\newcommand{\bd}{\mathbf{d}}

\title[]{Moments of Student's t-distribution: A Unified Approach}
\author[J.L. Kirkby]{J. Lars Kirkby }
\address{School of Industrial and Systems Engineering \\
Georgia Institute of Technology\\
Atlanta, GA 30318 \\
United States}
\email{jkirkby3@gatech.edu}

\author[D.H. Nguyen]{Dang H. Nguyen }
\address{Department of Mathematics \\
University of Alabama \\
Tuscaloosa, AL 35487\\
United States}
\email{dangnh.maths@gmail.com}

\author[D. Nguyen]{Duy Nguyen }
\address{Department of Mathematics\\
Marist College\\
Poughkeepsie, NY 12601\\
United States
}
\email{nducduy@gmail.com}


\begin{document}
\maketitle
\begin{abstract}
In this paper, we derive \emph{new} closed form formulae for moments of (generalized) Student's t-distribution 
in the \emph{one dimensional} case as well as in \emph{higher dimensions} through a unified probability framework. Interestingly,  the closed form expressions for the moments of Student's t-distribution can be written in terms of the familiar Gamma function, Kummer's confluent hypergeometric function, and the hypergeometric function.  This work aims to provide a concise and unified treatment of the moments for this important distribution.
\end{abstract}
\tableofcontents

\newpage
\section{Introduction}
In probability and statistics,
the location (e.g., mean), spread (e.g, standard deviation),
 skewness, and kurtosis play an important role in the modeling of random processes. One often
uses the mean and standard deviation to construct confidence
intervals or conduct hypothesis testing, and 
significant skewness or kurtosis of a data set indicates deviations from normality. Moreover, moment matching algorithms are among the most widely used fitting procedures in practice. As a result, it is important
to be able to find the moments of a given distribution. In his popular note, \cite{winkelbauer2012moments}
gave 
the closed form formulae
for the moments as well as absolute moments
of a normal distribution $N(\mu,\sigma^2)$.
The obtained results are beautiful
and have been well received. Recently, \cite{ogasawara2020unified}
provides a unified, non-recursive formulae for moments of normal distribution
with strip truncation. Also see \citep{nguyen2021probabilistic,skorski2024handy} for the binomial family.
Given the close relationship between the normal and Student's t-distributions, a natural question arises: Can we derive 
similar formulae for the family
of Student's t-distributions?
From the authors' best knowledge,
no such set of formulae exist
for (generalized) Student's t-distributions.
The purpose of this note
is to provide a complete set of closed form formulae
for raw moments, central raw moments, absolute moments, and central absolute moments
for (generalized) Student's t-distributions in the one-dimensional case and $n$-dimensional case.
In particular, the formulae given in  \eqref{musigmamoment} - \eqref{musigmaabscentralmoment} and Proposition \ref{proposition3} are new in the literature.
 In this sense, we unify existing results and provide extensions to higher dimensions, within a common probabilistic framework.\\

\noindent\textbf{Notation}: For later use,
we denote the probability density function (pdf) of a Gamma distribution  with parameters
$\alpha>0,\beta>0$ by
$$
\text{Gamma}(x|\alpha,\beta)
=\displaystyle\frac{\beta^{\alpha}}{\Gamma(\alpha)} x^{\alpha-1}e^{-\beta x},\quad x\in (0,\infty).
$$
Similarly,
the probability density function of a normal distribution
$X\sim N(\mu,\sigma^2)$ is denoted by
$$
N(x|\mu,\sigma^2)=\displaystyle\frac{1}{\sqrt{2\pi}\sigma}
\exp\left(-\displaystyle\frac{(x-\mu)^2}{2\sigma^2} \right), \quad x\in(-\infty,+\infty).
$$
This is extended naturally to higher dimensional cases.

We will also require two common special functions.
The Kummer's confluent hypergeometric function is defined by
$$
K(\alpha,\gamma; z)\equiv {}_1F_1(\alpha,\gamma; z)=\displaystyle\sum_{n=0}^{\infty}\frac{\alpha^{\overline{n}}z^n}{\gamma^{\overline{n}}n!}.
$$
The hypergeometric function is defined by
$$
{}_2F_1(a,b,c; z)=\sum_{n=0}^{\infty}\frac{a^{\bar n} b^{\bar n}}{c^{\bar n}}
\cdot\frac{z^n}{n!},
$$
where 
\begin{equation*}
a^{\bar n}=\frac{\Gamma(a+n)}{\Gamma(a)}=
\left\{
\begin{array}{ll}
1 & n=0,\\
a(a+1)\ldots(a+n-1) &n>0.
\end{array}
\right.
\end{equation*}

\section{Student's t-distribution: One dimensional case}
Recall the probability density function (pdf)
of a standard Student's t-distribution with {$\nu\in\{1,2,3,\ldots\}$ degrees of freedom}, denoted by  $ St(t|  0, 1,\nu)$, is given by
\begin{equation}
St(t|  0, 1,\nu)=\displaystyle\frac{\Gamma(\frac{\nu+1}{2})}{\Gamma(\frac{\nu}{2})}\frac{1}{\sqrt{\nu\pi}}
\left( 1+\frac{t^2}{\nu} \right)^{-\frac{\nu+1}{2}},\quad -\infty<t<\infty,
\end{equation}
where the Gamma function is defined as
$$\Gamma(z)=\displaystyle\int_0^{\infty} t^{z-1}e^{-t}dt.$$

\noindent {More generally, the probability density function of
a \textit{location-scale} (or generalized) Student's t-distribution with $\nu>0$ degrees of freedom is denoted by,}
\begin{equation}\label{t-genversion}
 St(t|  \mu, \sigma,\nu)=\displaystyle\frac{\Gamma(\frac{\nu+1}{2})}{\Gamma(\frac{\nu}{2})}
\left(\frac{\sigma}{{\nu\pi}}\right)^{\frac{1}{2}}
\left( 1+\frac{\sigma}{\nu}\left({t-\mu} \right)^2 \right)^{-\frac{\nu+1}{2}},\quad -\infty<t<\infty,
\end{equation}
where $\mu\in(-\infty,\infty)$ is the location, $\sigma>0$
determines the scale,
and $\nu\in\{1,2,3,\ldots\}$
{represents the number of the degrees of freedom}.
The thickness of its tails is determined by the degrees of freedom.
When $\nu=1$, the pdf in \eqref{t-genversion}
reduces
to the pdf of $\text{Cauchy}(\mu,\sigma)$, {while
 the pdf in \eqref{t-genversion}
converges to the pdf of
the normal $N(t|\mu,(1/\sqrt{\sigma})^2)$ as $\nu\to\infty$.}

While the tails of the normal distribution decay at  an exponential rate,
the Student's t-distribution is heavy-tailed, with a polynomial decay rate. Because of this, the Student's t-distribution has been widely adopted
in robust data analysis including (non) linear regression \citep{lange1989robust}, sample selection models \citep{marchenko2012heckman}, and
linear mixed effect models \citep{pinheiro2001efficient}. It is also among the most widely applied distributions for financial risk modeling, see \cite{QRM15}, \cite{Shaw06},\cite{kwon2020distribution}. The reader is
invited to refer to \cite{kotz2004multivariate} for more.

The mean and variance 
of a Student's t-distribution $T$
are well known and can be found in
closed form by using the properties
of the Gamma function. {Specifically for $\nu>2$, we have (see, for example, \cite{jackman2009bayesian}):
$$
\mathbb E(T)=\mu, \text{Var}(T)=\frac{1}{\sigma}\frac{\nu}{\nu-2}.
$$}

However, for higher order raw or central moments, the calculation quickly
becomes tedious. 



{
We note that one can use the fact that the Student's t-distribution can be written as  $T=X/\sqrt{Z/\nu}$ where $X\sim N(0,1), Z\sim \chi_\nu^2$, $X,Z$
are independent. From there, one can derive the probability density function of $T$.
We adopt the mixture approach which is surprisingly simple and will be very useful in later derivations. It provides a representation of a conditional Student's t-distribution in terms of a normal distribution, see, for example, \cite{bishop2006pattern}. More specifically, we have,}
\begin{lem}\label{mixtureRepre}
{Assume that for $\nu>0$,  $\Lambda\sim \text{Gamma}(\lambda|\nu/2,\nu/2)$.
Additional, given $\Lambda=\lambda$, assume further that $T|\lambda$ is a normal distribution with mean $\mu$ and variance $1/(\sigma\lambda)$.
Then $T$ is a $St(t|\mu,\sigma,\nu)$ Student's t-distribution.}
\end{lem}
\textbf{Proof:}
As the proof is very concise, we reproduce it here for the reader's convenience.
Let $f_T(t)$ be the probability density function of $T$.
We have
\begin{align*}
f_T(t)=&\displaystyle\int_0^{\infty}N(t|\mu,\frac{1}{\sigma\lambda})\text{Gamma}(\lambda|\frac{\nu}{2},\frac{\nu}{2})d\lambda \\
&=\displaystyle\int_0^{\infty}\frac{\sqrt{\sigma\lambda}}{\sqrt{2\pi}}e^{-\frac{\sigma\lambda}{2}(t-\mu)^2}\frac{\nu^{\nu/2}}{2^{\nu/2}\Gamma(\nu/2)}\lambda^{\nu/2-1} e^{-\frac{\nu}{2}\lambda}d\lambda\\
&=\frac{\sqrt{\sigma}}{\sqrt{2\pi}}\frac{\nu^{\nu/2}}{2^{\nu/2}\Gamma(\nu/2)}
\frac{\Gamma(\frac{\nu+1}{2})}{(\frac{\nu}{2}+\frac{\sigma}{2}(t-\mu)^2)^{\frac{\nu+1}{2}}}
\displaystyle\int_0^{\infty}\text{Gamma}(\lambda|\frac{\nu+1}{2},\frac{\nu}{2}+\frac{\sigma}{2}(t-\mu)^2))d\lambda\\
&=\frac{\sqrt{\sigma}}{\sqrt{2\pi}}\frac{\nu^{\nu/2}}{2^{\nu/2}\Gamma(\nu/2)}
\frac{\Gamma(\frac{\nu+1}{2})}{(\frac{\nu}{2}+\frac{\sigma}{2}(t-\mu)^2)^{\frac{\nu+1}{2}}}\\
&=St(t|\mu,\sigma,\nu).
\end{align*}
This completes the proof of the Lemma.
\qed
\\
\indent The following results are well known:
\begin{thm}\label{momentsNormal}
{
We have
\begin{enumerate}
\item If $X\sim N(0,\sigma^2)$ 
then
$$
\mathbb E(X^m)=
\left\{
\begin{array}{cl}
0, &\mbox{if} \quad m=2k+1, k\in\mathbb N\\
\displaystyle\frac{\sigma^m m!}{2^{m/2}(m/2)!}, &\mbox{if}\quad m=2k, k\in\mathbb N.
\end{array}
\right.
$$
\item If $X\sim\text{Gamma}(\alpha,\beta)$, then $\mathbb E(X^{\nu})=\frac{\beta^{-\nu} \Gamma(\nu+\alpha)}{\Gamma(\alpha)}$ for $-\alpha<\nu\in \mathbb R$.
\end{enumerate}
}
\end{thm}


With this and Lemma \ref{mixtureRepre} above, we are able to find
moments of Student's t-distribution. More specifically, we have the following comprehensive theorem in one dimension.
\begin{thm}
For $k\in \mathbb N_+$, $0< k < \nu$, the following results hold:
\begin{enumerate}
\item For $ T\sim St(t|  0, 1,\nu)$, the raw and absolute moments satisfy
\begin{equation}
\label{zeroonemoment}
\mathbb{E}(T^k) =
\left\{
\begin{array}{ll}
\frac{\Gamma(\frac{k+1}{2})}{\sqrt{\pi}} \cdot \frac{\nu^{k/2}}{\prod_{i=1}^{k/2} (\frac{\nu}{2}-i)}
,& k \text{ even},\\
0, & k \text{ odd}.
\end{array}
\right.
\end{equation}
\begin{equation}
\label{zerooneabsmoment}
\mathbb{E}(|T|^k)=
\frac{\nu^{k/2}\Gamma((k+1)/2)\Gamma((\nu-k)/2) }{\sqrt{\pi}\Gamma(\nu/2)}.
\end{equation}

\item If $ T\sim St(t|  \mu, \sigma,\nu)$, the raw moments satisfy
\begin{equation}
\label{musigmamoment}
\mathbb{E}(T^k)= \left\{
\begin{array}{ll}
(\nu/\sigma)^{k/2}
\displaystyle\frac{\Gamma(\frac{k+1}{2})}{\sqrt{\pi}}\frac{\Gamma(\frac{\nu}{2}-\frac{k}{2})}{\Gamma(\frac{\nu}{2})} {}_2F_1(-\frac{k}{2},\frac{\nu}{2}-\frac{k}{2},\frac{1}{2};-\frac{\mu^2\sigma}{\nu}),
 & k \text{ even},\\
2\mu(\nu/\sigma)^{(k-1)/2}\frac{\Gamma(\frac{k}{2}+1)}{\sqrt{\pi}} \frac{\Gamma(\frac{\nu}{2}-\frac{k-1}{2})}{\Gamma(\frac{\nu}{2})}  {}_2F_1(\frac{1-k}{2},\frac{\nu}{2}-\frac{k-1}{2},\frac{3}{2};-\frac{\mu^2\sigma}{\nu}), & k \text{ odd}.
\end{array}
\right.
\end{equation}
\begin{equation}
\label{musigmacentralmoment}
\mathbb E((T-\mu)^k)=\frac{(1+(-1)^k)}{2}(\nu/\sigma)^{k/2}\displaystyle\frac{\Gamma(\frac{k+1}{2})}{\sqrt{\pi}}\frac{\Gamma(\frac{\nu-k}{2})}{\Gamma(\frac{\nu}{2})}.
\end{equation}

\item If $ T\sim St(t|  \mu, \sigma,\nu)$, the absolute moments satisfy
\begin{equation}
\label{musigmaabsmoment}
\mathbb E(|T|^k)=(\nu/\sigma)^{k/2}
\displaystyle\frac{\Gamma(\frac{k+1}{2})}{\sqrt{\pi}}\frac{\Gamma(\frac{\nu}{2}-\frac{k}{2})}{\Gamma(\frac{\nu}{2})} {}_2F_1(-\frac{k}{2},\frac{\nu}{2}-\frac{k}{2},\frac{1}{2};-\frac{\mu^2\sigma}{\nu}).
\end{equation}
\begin{equation}
\label{musigmaabscentralmoment}
\mathbb E(|T-\mu|^k)=\displaystyle(\nu/\sigma)^{k/2}\frac{\Gamma(\frac{k+1}{2})}{\sqrt{\pi}}\frac{\Gamma(\frac{\nu-k}{2})}{\Gamma(\frac{\nu}{2})}.
\end{equation}
\end{enumerate}
In general, the moments are undefined when $k\geq \nu$.
\end{thm}
\noindent\textbf{Proof:}
First assume that $T\sim St(t|0,1,\nu)$; we will find 
$\mathbb{E}(|T|^k)$. The proof for $\mathbb{E}(T^k)$ follows from similar ideas in combination with the result obtained in Theorem \ref{momentsNormal}.
{Assume  $\Lambda\sim \text{Gamma}(\lambda|\nu/2,\nu/2)$.
Additional, given $\Lambda=\lambda$, assume further that $T|\lambda$ is a normal distribution with mean $0$ and variance $1/\lambda$.}
From the equation (17) in \cite{winkelbauer2012moments}, we have
\begin{align*}
\mathbb{E}(|T|^k|\lambda) &= \int \limits_{-\infty}^\infty |t|^k \text{ N}(t|0,\tfrac{1}{\lambda}) \text{ } dt =\frac{1}{\lambda^{k/2}}2^{k/2}
\displaystyle\frac{\Gamma(\frac{k+1}{2})}{\sqrt{\pi}}
K\left( -\frac{k}{2},\frac{1}{2};0 \right).
\end{align*}
Hence we have
\begin{align*}
&\mathbb{E}(|T|^k)=\mathbb E( \mathbb{E}(|T|^k|\lambda) )\\
&=\int \limits_{0}^\infty \frac{1}{\lambda^{k/2}}2^{k/2}
\displaystyle\frac{\Gamma(\frac{k+1}{2})}{\sqrt{\pi}}
K\left( -\frac{k}{2},\frac{1}{2};0 \right)
\cdot \frac{\nu^{\nu/2}}{2^{\nu/2} \Gamma(\tfrac{\nu}{2})} \lambda^{\nu/2-1} \exp \Big( -
\frac{\nu}{2} \lambda \Big) \text{ } d\lambda\\
&= 2^{k/2}
\displaystyle\frac{\Gamma(\frac{k+1}{2})}{\sqrt{\pi}}
K\left( -\frac{k}{2},\frac{1}{2};0 \right)
\cdot \frac{\nu^{\nu/2}}{2^{\nu/2} \Gamma(\tfrac{\nu}{2})} \int \limits_{0}^\infty \lambda^{\nu/2-1-k/2} \exp \Big( -
\frac{\nu}{2} \lambda \Big) \text{ } d\lambda\\
&= 2^{k/2}
\displaystyle\frac{\Gamma(\frac{k+1}{2})}{\sqrt{\pi}}
K\left( -\frac{k}{2},\frac{1}{2};0 \right)
\cdot \frac{\nu^{\nu/2}}{2^{\nu/2} \Gamma(\tfrac{\nu}{2})}
\cdot\frac{\Gamma(\frac{\nu-k}{2})}{(\frac{\nu}{2})^{\frac{\nu-k}{2}}}
 \int \limits_{0}^\infty \frac{(\frac{\nu}{2})^{\frac{\nu-k}{2}}}{\Gamma(\frac{\nu-k}{2})}\lambda^{(\nu-k)/2-1} \exp \Big( -
\frac{\nu}{2} \lambda \Big) \text{ } d\lambda\\
&=2^{k/2}
\displaystyle\frac{\Gamma(\frac{k+1}{2})}{\sqrt{\pi}}
K\left( -\frac{k}{2},\frac{1}{2};0 \right)
\cdot \frac{\nu^{\nu/2}}{2^{\nu/2} \Gamma(\tfrac{\nu}{2})}
\cdot\frac{\Gamma(\frac{\nu-k}{2})}{(\nu/2)^{\frac{\nu-k}{2}}}\\
&=\frac{\nu^{k/2}\Gamma((k+1)/2)\Gamma((\nu-k)/2) }{\sqrt{\pi}\Gamma(\nu/2)},
\end{align*}
where we have used the fact that $K\left( -\frac{k}{2},\frac{1}{2};0 \right)=1$. \\
\indent Next, {assume that $T\sim St(t|\mu,\sigma,\nu)$ and   $\Lambda\sim \text{Gamma}(\lambda|\nu/2,\nu/2)$.
Additional, given $\Lambda=\lambda$, assume further that $T|\lambda$ is a normal distribution with mean $\mu$ and variance $1/(\sigma\lambda)$.} Using the following facts (obtained in \cite{winkelbauer2012moments})
\begin{align*}
\mathbb{E}((T-\mu)^k|\lambda) &= \int \limits_{-\infty}^\infty (t-\mu)^k \text{ N}(t|\mu,\tfrac{1}{\lambda \sigma}) \text{ } dt =(1+(-1)^k)\frac{1}{\lambda^{k/2}}2^{k/2-1}\sigma^{-k/2}\displaystyle\frac{\Gamma(\frac{k+1}{2})}{\sqrt{\pi}},
\end{align*}
and
\begin{align*}
\mathbb{E}(|T-\mu|^k|\lambda) &= \int \limits_{-\infty}^\infty |t-\mu|^k \text{ N}(t|\mu,\tfrac{1}{\lambda \sigma}) \text{ } dt =\frac{\sigma^{-k/2}}{\lambda^{k/2}}2^{k/2}\displaystyle\frac{\Gamma(\frac{k+1}{2})}{\sqrt{\pi}},
\end{align*}
the derivations for $\mathbb{E}((T-\mu)^k)$ and $\mathbb{E}(|T-\mu|^k)$ follow similarly.\\
\indent {Next assume that $T\sim St(t|\mu,\sigma,\nu)$.  We would like to compute the absolute raw moment of $T$.} From the equation (17) in \cite{winkelbauer2012moments}, we have
\begin{align*}
\mathbb{E}(|T|^k|\lambda) &= \int \limits_{-\infty}^\infty |t|^k \text{ N}(t|\mu,\tfrac{1}{\lambda\sigma}) \text{ } dt =\frac{1}{\lambda^{k/2}}2^{k/2}\sigma^{-k/2}
\displaystyle\frac{\Gamma(\frac{k+1}{2})}{\sqrt{\pi}}
K\left( -\frac{k}{2},\frac{1}{2};-\frac{\mu^2}{2}\sigma\lambda \right).
\end{align*}
Hence, using Part 2) of Theorem \ref{momentsNormal}, we have {for $k<\nu$}
\begin{align*}
&\mathbb{E}(|T|^k)=\mathbb E( \mathbb{E}(|T|^k|\lambda) )\\
&=\int \frac{1}{\lambda^{k/2}}2^{k/2}\sigma^{-k/2}
\displaystyle\frac{\Gamma(\frac{k+1}{2})}{\sqrt{\pi}}
K\left( -\frac{k}{2},\frac{1}{2};-\frac{\mu^2}{2}\sigma\lambda \right)\text{Gamma}(\lambda|\frac{\nu}{2},\frac{\nu}{2})d\lambda\\
&=\sigma^{-k/2}2^{k/2}
\displaystyle\frac{\Gamma(\frac{k+1}{2})}{\sqrt{\pi}} 
\sum_{n=0}^{\infty}\frac{(-k/2)^{\bar n}}{(1/2)^{\bar n}}\frac{(-\mu^2/2)^n\sigma^n }{n!} \int \lambda^{n-k/2}\text{Gamma}(\lambda|\frac{\nu}{2},\frac{\nu}{2})d\lambda\\
&=\sigma^{-k/2}2^{k/2}
\displaystyle\frac{\Gamma(\frac{k+1}{2})}{\sqrt{\pi}} 
\sum_{n=0}^{\infty}\frac{(-k/2)^{\bar n}}{(1/2)^{\bar n}}\frac{(-\mu^2/2)^n\sigma^n }{n!} (\nu/2)^{-n+k/2}\frac{\Gamma(n-k/2+\nu/2)}{\Gamma(\nu/2)}\\
&=\sigma^{-k/2}2^{k/2}
\displaystyle\frac{\Gamma(\frac{k+1}{2})}{\sqrt{\pi}} \frac{\Gamma(\frac{\nu}{2}-\frac{k}{2})}{\Gamma(\frac{\nu}{2})}
\sum_{n=0}^{\infty}\frac{(-k/2)^{\bar n}}{(1/2)^{\bar n}}\frac{(-\mu^2/2)^n\sigma^n }{n!} (\nu/2)^{-n+k/2}(\frac{\nu}{2}-\frac{k}{2})^{\bar n}\\
&=\sigma^{-k/2}2^{k/2}(\nu/2)^{k/2}
\displaystyle\frac{\Gamma(\frac{k+1}{2})}{\sqrt{\pi}} \frac{\Gamma(\frac{\nu}{2}-\frac{k}{2})}{\Gamma(\frac{\nu}{2})}
\sum_{n=0}^{\infty}\frac{(-k/2)^{\bar n}}{(1/2)^{\bar n}}\frac{(-\mu^2/2)^n\sigma^n }{n!} (\nu/2)^{-n}(\frac{\nu}{2}-\frac{k}{2})^{\bar n}\\
&= (\nu/\sigma)^{k/2}
\displaystyle\frac{\Gamma(\frac{k+1}{2})}{\sqrt{\pi}}\frac{\Gamma(\frac{\nu}{2}-\frac{k}{2})}{\Gamma(\frac{\nu}{2})} {}_2F_1(-\frac{k}{2},\frac{\nu}{2}-\frac{k}{2},\frac{1}{2};-\frac{\mu^2\sigma}{\nu}).
\end{align*}
Lastly, from the equation (12) in \cite{winkelbauer2012moments}, we have
\begin{align*}
\mathbb{E}(T^k|\lambda)&= \int \limits_{-\infty}^\infty t^k \text{ N}(t|\mu,\tfrac{1}{\lambda \sigma}) \text{ } dt\\
&=\left\{
\begin{array}{ll}
\sigma^{-k/2}2^{k/2}\frac{\Gamma(\frac{k+1}{2})}{\sqrt{\pi}}\frac{1}{\lambda^{k/2}}K(-\frac{k}{2},\frac{1}{2};
-\frac{\mu^2}{2}\sigma\lambda) & k \text{ even},\\
\mu\sigma^{-(k-1)/2}2^{(k+1)/2}\frac{\Gamma(\frac{k}{2}+1)}{\sqrt{\pi}}\frac{1}{\lambda^{(k-1)/2}}K(\frac{1-k}{2},\frac{3}{2};
-\frac{\mu^2}{2}\sigma\lambda) &\text{ } k \text{ odd}.
\end{array}
\right.
\end{align*}
Similar to the calculations done for $\mathbb{E}(|T|^k)$, we have
\begin{equation*}
\mathbb{E}(T^k)= \left\{
\begin{array}{ll}
(\nu/\sigma)^{k/2}
\displaystyle\frac{\Gamma(\frac{k+1}{2})}{\sqrt{\pi}}\frac{\Gamma(\frac{\nu}{2}-\frac{k}{2})}{\Gamma(\frac{\nu}{2})} {}_2F_1(-\frac{k}{2},\frac{\nu}{2}-\frac{k}{2},\frac{1}{2};-\frac{\mu^2\sigma}{\nu}).
 & k \text{ even},\\
2\mu(\nu/\sigma)^{(k-1)/2}\frac{\Gamma(\frac{k}{2}+1)}{\sqrt{\pi}} \frac{\Gamma(\frac{\nu}{2}-\frac{k-1}{2})}{\Gamma(\frac{\nu}{2})}  {}_2F_1(\frac{1-k}{2},\frac{\nu}{2}-\frac{k-1}{2},\frac{3}{2};-\frac{\mu^2\sigma}{\nu}) & k \text{ odd}.
\end{array}
\right.
\end{equation*}
This completes the proof of the theorem.
\qed
\begin{rmk}
\text{}
\begin{enumerate}
\item The formulae given in  \eqref{musigmamoment} - \eqref{musigmaabscentralmoment} are new in the literature. Also when $ T\sim St(t|  0, 1,\nu)$, $\mathbb E(T^k)$ is well known. Moreover, one can directly use
the definition to find $\mathbb E(|T|^k)$
through the class of $\beta$-functions defined in Section 6.2 of \cite{abramowitz1948handbook} and arrive
at the same formula.
However, this direct approach no longer works for expectations of the form
$\mathbb E(|T|^{k})$, $\mathbb E(T^{k})$ when $ T\sim St(t|  \mu, \sigma,\nu)$, or for higher dimensional moments
considered in Section \ref{HigherDimcase}. {Also, clearly \eqref{musigmamoment} is reduced to \eqref{zeroonemoment} and \eqref{musigmaabsmoment} is reduced to \eqref{zerooneabsmoment} when $\mu=0$ and $\sigma=1$.} 
\item If $ T\sim St(t|  \mu, \sigma,\nu)$, and once $\mathbb E((T-\mu)^i), 0\leq i\leq k$ have been computed, we can use them to compute
$\mathbb E(T^k)$ for $k<\nu$ using the expansion
$$
\mathbb E(T^k)=\mathbb E((T-\mu+\mu)^k)=\sum_{i=0}^{k}\mu^{k-i}{k\choose i}\mathbb E((T-\mu)^i).
$$
\item {We also note that an alternative proof of the central moment formula \eqref{musigmacentralmoment} was later provided in \cite{bignozzi2024inter} based on a recursive formula.
}
\end{enumerate}
\end{rmk}

\section{Higher dimensional case}\label{HigherDimcase}
Now we consider the case when $n\geq 2$. Denote $\bm t=(t_1,t_2,\ldots,t_n)\in\mathbb R^{n}$.
Denote the pdf of $n$-dimensional Normal random variable as
 $$
N(\bm x|\bm\mu, \bm\Sigma)
=\frac{1}{(2\pi)^{n/2}|\bm\Sigma|^{\frac{1}{2}}}
e^{-\frac{1}{2}(\bm x-\bm \mu)^T \bm\Sigma^{-1}(\bm x-\bm \mu)},\bm x\in\mathbb R^n,
$$
where $\bm\mu\in\mathbb R^n$ and $|\bm\Sigma|$ is the determinant of {the $n\times n$ - symmetric positive definite matrix} $\bm\Sigma$. Similar the 1-dimensional case,  we have
the probability density of $n$-dimensional Student's t-distribution is defined as,
\begin{equation}\label{densitythroughnormal}
St(\bm t| \bm \mu, \bm \Sigma,\nu)
=\displaystyle\int_0^{\infty} N(\bm t|\bm \mu, (\lambda\bm\Sigma)^{-1})\text{Gamma}(\lambda|\frac{\nu}{2}
, \frac{\nu}{2})d\lambda,
\end{equation}
where  $\bm\mu$ is called the location, $\bm\Sigma$ is the  scale matrix, and $\nu$
is the degrees of freedom parameter. Similar to Lemma \ref{mixtureRepre}, we have 

\begin{align*}
    &St(\bm t \mid {\bm\mu}, \boldsymbol{\Sigma}, \nu)  = \int_0^\infty N(\bm t \mid {\bm\mu}, (\lambda \boldsymbol{\Sigma})^{-1})  \text{Gamma}\left(\lambda \mid \frac{\nu}{2}, \frac{\nu}{2}\right) d\lambda \\
    & = \int_0^\infty \frac{|\lambda \boldsymbol{\Sigma}|^{1/2}}{(2\pi)^{n/2}}  \exp \left\{-\frac{1}{2}(\bm t-\bm\mu)^T(\lambda \boldsymbol{\Sigma})(\bm t-\bm\mu) - \frac{\nu \lambda}{2}\right\} 
     \frac{1}{\Gamma(\nu/2)} \left(\frac{\nu}{2}\right)^{\nu/2} \lambda^{\nu/2-1} d\lambda \\
    & = \frac{(\nu/2)^{\nu/2} |\boldsymbol{\Sigma}|^{1/2}}{(2\pi)^{n/2} \Gamma(\nu/2)} \int_0^\infty \exp \left\{-\frac{1}{2}(\bm t-\bm\mu)^T(\lambda \boldsymbol{\Sigma})(\bm t-\bm\mu) - \frac{\nu \lambda}{2}\right\}  \lambda^{n/2+\nu/2-1} d\lambda.
\end{align*}
Let's define 
\begin{align*}
    \Delta^2 & = (\bm t-\bm\mu)^T \boldsymbol{\Sigma} (\bm t-\bm\mu), \\
    z & = \frac{\lambda}{2}(\Delta^2 + \nu),
\end{align*}
\noindent then we have

\begin{align}
\label{n-t-dist-pdf}
    St(x \mid {\bm\mu}, \boldsymbol{\Sigma}, \nu) & = \frac{(\nu/2)^{\nu/2} |\boldsymbol{\Sigma}|^{1/2}}{(2\pi)^{n/2} \Gamma(\nu/2)} \int_0^\infty \exp(-z) \left(\frac{2z}{\Delta^2 + \nu}\right)^{n/2+\nu/2-1} \cdot \frac{2}{\Delta^2 + \nu} dz\notag \\
    & = \frac{(\nu/2)^{\nu/2} |\boldsymbol{\Sigma}|^{1/2}}{(2\pi)^{n/2} \Gamma(\nu/2)} \left(\frac{2}{\Delta^2 + \nu}\right)^{n/2+\nu/2} \int_0^\infty \exp(-z) z^{n/2+\nu/2-1} dz \notag\\
    & = \frac{(\nu/2)^{\nu/2} |\boldsymbol{\Sigma}|^{1/2}}{(2\pi)^{n/2} \Gamma(\nu/2)} \left(\frac{2}{\Delta^2 + \nu}\right)^{n/2+\nu/2}\Gamma(\frac{\nu+n}{2})\notag\\
    &=\displaystyle
\frac{\Gamma(\frac{\nu+n}{2})}{\Gamma(\frac{\nu}{2})}
\frac{|\bm\Sigma|^{\frac{1}{2}}}{(\nu\pi)^{\frac{n}{2}}}
\left(1+\frac{1}{\nu}(\bm t-\bm\mu)^T\bm\Sigma(\bm t -\bm\mu) \right)^{-\frac{\nu+n}{2}}.
\end{align}


\noindent Note that in the standardized case of
$\bm \mu=\bm 0$ and $\bm\Sigma =\bm I$, the representation in \eqref{densitythroughnormal} is reduced to 
\begin{equation}\label{densitythroughnormal1}
St(\bm t| \bm 0, \bm I,\nu)
=\displaystyle\int_0^{\infty} N(\bm t|\bm 0, \frac{1}{\lambda}\bm I)\text{Gamma}(\lambda|\frac{\nu}{2}
, \frac{\nu}{2})d\lambda.
\end{equation}


\noindent Let's write $\bm T=(T_1,T_2,\ldots,T_n)$, and  $\bm k=(k_1,k_2,\ldots,k_n)$ with $0\leq k_i\in\mathbb N$.
{For $\nu>2$, it is known that (see, for example, \cite{bishop2006pattern}),
$$
\mathbb E[\bm T]=\bm\mu, \text{Cov}(\bm T)=\frac{\nu}{\nu-2}\bm\Sigma^{-1}.
$$
}
For the rest of this section, we are interested in higher moments of $\bm T$. 
The $\bm k$  moment of $\bm T$
is defined as
$$
\mathbb E(\bm T^{\bm k})=\bm\int t_1^{k_1}t_2^{k_2}\ldots t_n^{k_n} \cdot St(\bm t| \bm \mu, \bm \Sigma,\nu)dt_1\ldots dt_n.
$$
Similarly,
$$
\mathbb E( |\bm T|^{\bm k})=\bm\int |t_1|^{k_1}|t_2|^{k_2}\ldots |t_n|^{k_n}\cdot St(\bm t| \bm \mu, \bm \Sigma,\nu)dt_1\ldots dt_n.
$$
{To simplify the notations,
in the following, we use $\sum k_i,\prod k_i$
to denote $\sum_{i=1}^n k_i,\prod_{i=1}^n k_i$, respectively}. From the authors' best knowledge, the following are new:
\begin{thm}
\label{proposition3}
For $\sum k_i< \nu$, we have 
\begin{enumerate}
\item If $\bm T\sim St(\bm t| \bm 0, \bm I,\nu)$ then
\begin{itemize}
\item The raw moments satisfy
$$
\mathbb E(\bm T^{\bm k})=\left\{
\begin{array}{ll}
0, &\mbox{if} \quad \text{at least one} \text{ } k_i\text{ is odd},\\
\displaystyle \displaystyle\nu^{\frac{\sum k_i}{2}}\frac{ \Gamma(\frac{\nu-\sum k_i}{2})}{ 
\Gamma(\frac{\nu}{2})} \frac{ \prod (k_i)!}{2^{(\sum k_i)}\prod(k_i/2)!}, &\mbox{if} \text{ all } k_i\text{ are even}.
\end{array}
\right.
$$
\item  The absolute moments satisfy
$$\mathbb E(|\bm T|^{\bm k})=\displaystyle\nu^{\frac{\sum k_i}{2}}\frac{ \Gamma(\frac{\nu-\sum k_i}{2})}{ 
\Gamma(\frac{\nu}{2})}\prod\frac{\Gamma(\frac{k_i+1}{2})}{\sqrt{\pi}}.$$
\end{itemize}
\item If $\bm T\sim St(\bm t| \bm \mu, \bm \Sigma,\nu)$, denote  $\bm\Sigma^{-1}=(\overline{\sigma}_{ij})$ and $\bm e_i=(0,\ldots,1,\ldots,0)$ - the $i$th unit vector of $\mathbb R^n$. Then we have the following recursive
formula to compute the moments of $\bm T$:
\begin{align*}
\mathbb E(\bm T^{\bm k+\bm e_i})
=\mu_i\mathbb E(\bm T^{\bm k})
+\frac{\nu}{\nu-2}\sum_{j=1}^{n}\overline{\sigma}_{ij}k_j\mathbb E(\bm T^{\bm k-\bm e_j}).
\end{align*}

\end{enumerate}

\end{thm}
\noindent\textbf{Proof}:\\
For 1), first from \eqref{densitythroughnormal1}, we have
\begin{align*}
\mathbb E(\bm T^{\bm k})=\int_0^{\infty}\mathbb E(\bm X^{\bm k}|\bm 0,\frac{1}{t}\bm I)\text{Gamma}(t|\frac{\nu}{2},\frac{\nu}{2})dt,
\end{align*}
where $\E(\bm X^{\bm k}|\bm 0,\frac{1}{t}\bm I)$ is the $\bm k$ moment  of a $N(\bm 0,\frac{1}{t}\bm I)$.
Using Theorem \ref{momentsNormal}, we have
$$
\E(\bm X^{\bm k}|\bm 0,\frac{1}{t}\bm I)=
\prod_{i=1}^{n} \E( X_i^{k_i}|0,\frac{1}{t})=
\left\{
\begin{array}{cl}
0, &\mbox{if} \quad \text{at least one }  k_i\text{ is odd},\\
\displaystyle\frac{t^{-\sum k_i/2} \prod (k_i)!}{2^{(\sum k_i)/2}\prod(k_i/2)!}, &\mbox{if} \quad \text{all } k_i\text{ are even}.
\end{array}
\right.
$$ 
As a result,
\begin{align*}
\mathbb E(\bm T^{\bm k})&=
\left\{
\begin{array}{ll}
0, &\mbox{if} \text{ at least one }  k_i\text{ is odd},\\
\displaystyle\frac{ \prod (k_i)!}{2^{(\sum k_i)/2}\prod(k_i/2)!}\int_0^{\infty} t^{-\sum k_i/2}\text{Gamma}(t|\frac{\nu}{2},\frac{\nu}{2})dt, &\mbox{if} \text{ all } k_i\text{ are even}.
\end{array}
\right.\\
&=\left\{
\begin{array}{ll}
0, &\mbox{if}  \text{ at least  one }  k_i \text{ is odd},\\
\displaystyle \displaystyle\nu^{\frac{\sum k_i}{2}}\frac{ \Gamma(\frac{\nu-\sum k_i}{2})}{ 
\Gamma(\frac{\nu}{2})} \frac{ \prod (k_i)!}{2^{(\sum k_i)}\prod(k_i/2)!}, &\mbox{if} \text{ all } k_i\text{ are even}.
\end{array}
\right.
\end{align*}

\noindent Similarly, we have
\begin{align*}
\mathbb E(|\bm T^{\bm k}|)=\int_0^{\infty}\mathbb E(|X_1|^{k_1}|X_2|^{k_2}\ldots |X_n|^{k_n}|\bm 0,\frac{1}{t}\bm I)\text{Gamma}(t|\frac{\nu}{2},\frac{\nu}{2})dt
\end{align*}
where
\begin{align*}
\mathbb E(|X_1|^{k_1}|X_2|^{k_2}\ldots |X_n|^{k_n}|\bm 0,\frac{1}{t}\bm I)=
\prod_{i=1}^{n} \E( |X_i|^{k_i}|0,\frac{1}{t})=\prod \frac{1}{t^{k_i/2}}2^{k_i/2}
\displaystyle\frac{\Gamma(\frac{k_i+1}{2})}{\sqrt{\pi}}.
\end{align*}
Therefore,
\begin{align*}
\mathbb E( |\bm T^{\bm k}|)&=2^{\sum k_i/2}
\displaystyle\prod\frac{\Gamma(\frac{k_i+1}{2})}{\sqrt{\pi}}\int_0^{\infty}t^{-\sum k_i/2}\text{Gamma}(t|\frac{\nu}{2},\frac{\nu}{2})dt
\\
&=\displaystyle\nu^{\frac{\sum k_i}{2}}\frac{ \Gamma(\frac{\nu-\sum k_i}{2})}{ 
\Gamma(\frac{\nu}{2})}\prod\frac{\Gamma(\frac{k_i+1}{2})}{\sqrt{\pi}}\quad \text{if}\quad \sum k_i< \nu.
\end{align*}
For 2), from \eqref{densitythroughnormal}
\begin{align}\label{recursiveEq}
\mathbb E(\bm T^{\bm k})=\int_0^{\infty}\mathbb E(\bm X^{\bm k}|\bm\mu,\frac{1}{t}\bm \Sigma^{-1})\text{Gamma}(t|\frac{\nu}{2},\frac{\nu}{2})dt,
\end{align}
where $\mathbb E(\bm X^{\bm k})\equiv \mathbb E(\bm X^{\bm k}|\bm\mu,\frac{1}{t}\bm \Sigma^{-1})$
is the $\bm k$ moment of $N(\bm\mu,\frac{1}{t}\bm \Sigma^{-1})$.
Recall the pdf of $N(\bm \mu,\frac{1}{t}\Sigma^{-1})$ is given by $$
N(\bm x|\bm\mu, \frac{1}{t}\bm\Sigma^{-1})
=\frac{1}{(2\pi)^{n/2}|\frac{1}{t}\bm\Sigma^{-1}|^{\frac{1}{2}}}
e^{-\frac{1}{2}(\bm x-\bm \mu)^T t\bm\Sigma(\bm x-\bm \mu)}.
$$
Similar to Theorem 1 in \cite{kan2017moments}, we have
\begin{align*}
-\frac{\partial N(\bm x|\bm\mu, \frac{1}{t}\bm\Sigma^{-1})}{\partial\bm x}
=t\bm\Sigma(\bm x-\bm \mu)N(\bm x|\bm\mu, \frac{1}{t}\bm\Sigma^{-1}).
\end{align*}
Hence
\begin{align*}
-\bm\int\bm x^{\bm k}\frac{\partial N(\bm x|\bm\mu, \frac{1}{t}\bm\Sigma^{-1})}{\partial\bm x}d\bm x
=\bm\int \bm x^{\bm k} t\bm\Sigma(\bm x-\bm \mu)N(\bm x|\bm\mu, \frac{1}{t}\bm\Sigma^{-1})d\bm x.
\end{align*}
By integration by parts, we arrive at
\begin{align*}
\bm\int k_j\bm x^{\bm k-\bm e_j}N(\bm x|\bm\mu, \frac{1}{t}\bm\Sigma^{-1})d\bm x
=\bm\int \bm x^{\bm k} t\bm\Sigma(\bm x-\bm \mu)N(\bm x|\bm\mu, \frac{1}{t}\bm\Sigma^{-1})d\bm x.
\end{align*}
Or equivalently
\begin{align*}
\bm\int \bm x^{\bm k} (\bm x-\bm \mu)N(\bm x|\bm\mu, \frac{1}{t}\bm\Sigma^{-1})d\bm x=
\frac{1}{t}\bm\Sigma^{-1}\bm\int k_j\bm x^{\bm k-\bm e_j}N(\bm x|\bm\mu, \frac{1}{t}\bm\Sigma^{-1})d\bm x.
\end{align*}
This in turn implies that
$$
\mathbb E(\bm X^{\bm k+\bm e_i})
=\mu_i\mathbb E(\bm X^{\bm k})
+\frac{1}{t}\sum_{j=1}^{n}\overline{\sigma}_{ij}k_j\mathbb E(\bm X^{\bm k-\bm e_j}).
$$
Plugging this into the equation \eqref{recursiveEq}, we have the following recursive equation
\begin{align*}
\mathbb E(\bm T^{\bm k+\bm e_i})
&=\mu_i\mathbb E(\bm T^{\bm k})
+\sum_{j=1}^{n}\bar{\sigma}_{ij}k_j\mathbb E(\bm T^{\bm k-\bm e_j})\int_0^{\infty}\frac{1}{t}\text{Gamma}(t|\frac{\nu}{2},\frac{\nu}{2})dt\\
&=\mu_i\mathbb E(\bm T^{\bm k})
+\frac{\nu}{2}\frac{\Gamma(\frac{\nu}{2}-1)}{\Gamma(\frac{\nu}{2})}
\sum_{j=1}^{n}\overline{\sigma}_{ij}k_j\mathbb E(\bm T^{\bm k-\bm e_j})\\
&=\mu_i\mathbb E(\bm T^{\bm k})
+\frac{\nu}{\nu-2}
\sum_{j=1}^{n}\overline{\sigma}_{ij}k_j\mathbb E(\bm T^{\bm k-\bm e_j}).
\end{align*}
This completes the proof of the theorem.
\qed

\vspace{0.5cm}

\indent Lastly, for $\bm a=(a_1,a_2,\ldots,a_n)$ and $\bm b=(b_1,b_2,\ldots,b_n)\in\mathbb R^n$, let
$\bm a_{(j)}$ is the vector obtained from $\bm a$ by deleting the $j$th element of $\bm a$.
For $\bm\Sigma=(\sigma_{ij})$, let $\sigma_i^2=\sigma_{ii}$ and  $\bm\Sigma_{i,(j)}$
stand for the $i$th row of $\bm\Sigma$ with its $j$th element removed.
Analogously, let $\bm\Sigma_{(i),(j)}$
stand for the matrix $\bm\Sigma$ with $i$th row and $j$th column removed.
\\
\indent Consider the following truncated 
$\bm k$ moment
$$
F^n_{\bm k}(\ba,\bb; \bm\mu,\bm\Sigma,\nu)
=\int_{\bm a}^{\bm b} \bm t^{\bm k} St(\bm x| \bm \mu, \bm \Sigma,\nu)d\bm t
\equiv\int_{a_1}^{ b_1}\ldots \int_{a_n}^{ b_n}  t_1^{ k_1}\ldots t_{n}^{k_n} St(\bm x| \bm \mu, \bm \Sigma,\nu)dt_1\ldots dt_n.
$$
We have
	\begin{equation}\label{7e1}
	\begin{aligned}
	F^n_{\bm k}(\ba,\bb; \bm\mu,\bm\Sigma,\nu) =\int_0^{\infty}\mathbb E\left[\1_{\{\ba\leq \BX\leq\bb\}}\bm X^{\bm k}\Big|\bm\mu,\frac{1}{t}\bm \Sigma^{-1}\right]\text{Gamma}(t|\frac{\nu}{2},\frac{\nu}{2})dt,
	\end{aligned}
	\end{equation}
where $\bm X^{\bm k}\sim N(\bm x|\bm\mu,\frac{1}{t}\bm \Sigma^{-1})$.
Using Theorem 1 in \cite{kan2017moments}, we have for $n>1$
\begin{align}
\E (X_{\bm k+\bm e_i}^{n};\bm a, \bm b,\bm \mu,\frac{1}{t}\bm\Sigma^{-1})&:=\E\left[\1_{\{\ba\leq \BX\leq\bb\}}\bm X^n_{\bm k+\bee_i}\Big|\bm\mu,\frac{1}{t}\bm \Sigma^{-1}\right]\nonumber\\
&=\mu_i\E\left[\1_{\{\ba\leq \BX\leq\bb\}}\bm X^n_{\bm k}\Big|\bm\mu,\frac{1}{t}\bm \Sigma^{-1}\right]
+\frac1t \bee_i^\top\BSigma^{-1}\bc_\bk,
\end{align}

\noindent where $\bc_\bk$ satisfies
\begin{align*}
\bc_{\bk,j}
&=k_j\E (X_{\bm k-\bm e_i}^{n};\bm a, \bm b,\bm \mu,\frac{1}{t}\bm\Sigma^{-1})
+a_j^{k_j}N(a_j|\mu_j,\frac{1}{t}\bar\sigma_j^2)
\E (X_{\bm k_{(j)}}^{n-1};\bm a_{(j)}, \bm b_{(j)},\widehat{\bm \mu}_j^{\bm a},\frac{1}{t}\widehat{\bm\Sigma}_j)\\
&-b_j^{k_j}N(b_j|\mu_j,\frac{1}{t}\bar\sigma_j^2)
\E (X_{\bm k_{(j)}}^{n-1};\bm a_{(j)}, \bm b_{(j)},\widehat{\bm \mu}_j^{\bm b},\frac{1}{t}\widehat{\bm\Sigma}_j), j=1,2,\ldots,n,
\end{align*}
with
\begin{equation}
\left\{
\begin{array}{l}
\widehat{\bm \mu}_j^{\bm a}=\bm \mu_{(j)}+\bm\Sigma^{-1}_{(j),j}\frac{a_j-\mu_j}{\overline\sigma^2_j},\\
\widehat{\bm \mu}_j^{\bm b}=\bm \mu_{(j)}+\bm\Sigma^{-1}_{(j),j}\frac{b_j-\mu_j}{\overline\sigma^2_j},\\
\widehat{\bm\Sigma}_j=\bm\Sigma^{-1}_{(j),(j)}-\frac{1}{\overline\sigma^2_j}\bm\Sigma^{-1}_{(j),j}\bm \Sigma^{-1}_{j,(j)}.
\end{array}
\right.
\end{equation}
Thus, we have the following recursive formula

$$
	F^n_{\bm k+\bee_i}(\ba,\bb; \bm\mu,\bm\Sigma,\nu)
	=\mu_i 	F^n_{\bm k}(\ba,\bb; \bm\mu,\bm\Sigma,\nu)
	+\frac{\nu}{\nu-2} \bee_i^\top\BSigma\bd_\bk,
$$
where
\begin{align*}
\bd_{\bk,j}=&k_j F^n_{\bm k-\bee_i}(\ba,\bb; \bm\mu,\bm\Sigma,\nu)
+a_j^{k_j} St(a_j|  \mu_j, \bar\sigma_j^2,\nu) 
F^{n-1}_{\bm k_{(j)}}(\ba_{(j)},\bb_{(j)}; \widehat\bmu^\ba_j,\widehat{\bm\Sigma},\nu)\\
&-b_j^{k_j} St(b_j|  \mu_j, \bar\sigma_j^2,\nu)
F^{n-1}_{\bm k_{(j)}}(\ba_{(j)},\bb_{(j)}; \widehat\bmu^{ \bb}_j,\widehat{\bm\Sigma},\nu), j=1,2,\ldots,n.
\end{align*}
Note by convention that the first term, second term, and third term in the expression of $\bd_{\bk,j}$ equal 0 when $k_j=0$, $a_j=\infty$, $b_j=-\infty$, respectively.

\section{Conclusion}
We derive
the closed form formulae
for the raw moments, absolute moments, and central moments
of Student's t-distribution with arbitrary degrees of freedom. We provide results in one and $n$-dimensions, which unify and extend the existing literature for the Student's t-distribution.
It would be interesting to investigate tail quantile approximations or asymptotic tail properties
of higher (generalized) {Student's t-distribution} as done in \cite{schluter2012tail} and \cite{finner2008asymptotic}.
We leave this as an interesting project for future studies.

\bibliographystyle{agsm}
\bibliography{cliquet}
\appendix

\end{document}